\documentclass{amsart}

\newcommand{\remark}{\noindent {\bf Remark. }}

\newcommand{\ws}{\hspace{4pt}}
\newtheorem{theorem}{Theorem}

\newtheorem{proposition}{Proposition}
\newtheorem{cor}{Corollary}

\usepackage[]{amssymb,amsopn}
\usepackage[]{amsmath}
\usepackage[]{eucal}
\usepackage[]{latexsym}

\makeatletter
\@namedef{subjclassname@2020}{%
  \textup{2020} Mathematics Subject Classification}
\makeatother

\begin{document}

\title[Compactness criteria]{Compactness criteria via Laguerre and Hankel transformations}
\author{\'A. P. Horv\'ath }

\thanks{Supported by the NKFIH-OTKA Grants K128922 and K132097.}

\begin{abstract}The aim of this paper is to prove Kolmogorov-Riesz type theorems via Bessel and Laguerre translations, and Pego-type theorems by the corresponding transformations.

\end{abstract}
\maketitle

\section{Introduction}

The populous family of compactness theorems is established by Arzel\'a and Ascoli. Our starting point, the Kolmogorov-Riesz theorem characterizes precompact (totally bounded) sets of $L^p(\mathbb{R}^n)$, see e.g. \cite{a}. Besides, such a characterization is interesting in itself, it has several applications to differential and integral equations. Compactness criteria were studied in particular non-standard function spaces, e.g. in Sobolev spaces \cite{hoh}, variable Lebesgue spaces \cite{gma}, or weighted variable exponent amalgam and Sobolev spaces \cite{au}, and also in more general circumstances, see e.g. \cite{f}, \cite{r} and \cite{gr}. Sudakov-type improvements of the classical theorem are derived in \cite{hohm} (see \cite{s} as well).\\
The useful criterion of compactness given by Pego via Fourier transformation (see \cite{pe}) has also great influence, see e.g. \cite{dfg}, \cite{g}, \cite{gk}. In the inspirating paper \cite{k}, compactness criterion is given by Laplace transformation.

Below we offer a new aspect of characterization of precompact sets in certain Banach spaces. Instead of deriving similar theorems by related, for instance  Mellin or cosine transformations, we investigate and extend the notion of equicontinuity. Motivated by the effect of translation, $\tau_yf(x)=f(x-y)$, on exponential functions, different translation operators were introduced by orthogonal systems $\{\varphi_n\}_n$ as $T_y\varphi_n(x)=\varphi_n(x)\varphi_n(y)$. Hereinafter we deal with Laguerre and Bessel translations.

In the next section a Riesz-Kolmogorov type theorem is derived by Laguerre translation to weighted $L^p$ spaces on the half-line. The corresponding transformation, the discrete Laguerre-Fourier transformation implies a remarkable simplification, since the structure of the corresponding $l^{p'}$ spaces are simpler than the original $L^p$ ones. In this section we also introduce and study Laguerre translation on sequences.\\
In the third section the method is presented by Bessel translation. The corresponding transformation is the Hankel (Bessel-Fourier) transformation which establishes a Pego-type theorem.

\medskip

\section{Laguerre translation method}

Laguerre translation is developed by product formulae, see e.g. \cite{w} and  by investigation of the related Cauchy problem, see e.g. \cite{bs}. We mention that it is a natural idea to handle Bessel and Laguerre translations in parallel, since the translated functions in both cases are the solutions to very similar Cauchy problems, $u_{xx}-u_{tt}+\frac{2\alpha+1}{x}u_x - \frac{2\alpha+1}{t}u_t-ru=0$, where $r(x,t)=0$ in Bessel case and $r(x,t)=x^2-t^2$ in Laguerre case. The derived convolution structure is examined e.g. in \cite{gm}. The norm of the translation operator ensures a maximum principle to the corresponding hyperbolic problem, and it implies Nikol'skii-type norm-estimations on the half-line, see \cite{adh}. Convolution is applicable to study of best approximation in certain spaces, see \cite{gm1}. In this section compactness is investigated by Laguerre translation. The main result of the section is the characterization of precompact sets in $L^2_\alpha$.

\subsection{Laguerre translation and precompact sets in  $L^p_\alpha$}

We introduce weighted $L^p$-spaces on the half-line similar to the ones given in the previous section.\\ Let $1\le p \le \infty$ and $\alpha>-1$.
$f: \mathbb{R}_+ \to \mathbb{R} \in L^p_\alpha$ if
$$\|f\|_{p,\alpha} := \left(\int_0^\infty \left|f(x) e^{-\frac{x}{2}}\right|^px^\alpha dx\right)^{\frac{1}{p}}=\left(\int_0^\infty |\tilde{f}|^pd\mu_\alpha\right)^{\frac{1}{p}}<\infty.$$
$$\|f\|_{\infty,\alpha}=\|\tilde{f}\|_\infty$$
independently of $\alpha$.

The translation operator acting on the space given above can be defined as follows, see e.g. \cite{gm} .\\
Let $\alpha>-\frac{1}{2}$
\begin{equation}\label{transl}T_t^\alpha(f,x)\end{equation}$$=c_{\alpha}\int_0^\pi f(x+t+2\sqrt{xt}\cos\vartheta)e^{-\sqrt{xt}\cos\vartheta}j_{\alpha-\frac{1}{2}}(\sqrt{xt}\sin\vartheta)\sin^{2\alpha}\vartheta d\vartheta,$$
where $c_{\alpha}=\frac{\Gamma(\alpha+1)}{\Gamma\left(\alpha+\frac{1}{2}\right)\Gamma\left(\frac{1}{2}\right)}$ and $j_{\alpha-\frac{1}{2}}$ is the entire Bessel function, cf. \eqref{B}.

By symmetry of the definition in $t,x \ge 0$,
$$T_t^\alpha(f,x)=T_x^\alpha(f,t),$$
and again by definition
$$T_t^\alpha(f,0)=f(t), \ws \ws \ws T_0^\alpha(f,x)=f(x).$$

Denoting by $L_n^\alpha(x)$ the $n^{th}$ Laguerre polynomial orthogonal on $(0,\infty)$ with respect to the Laguerre weight, that is
$$\int_0^\infty L_n^\alpha(x)L_k^\alpha(x)e^{-x}x^\alpha dx =\Gamma(\alpha+1)\binom{n+\alpha}{n}\delta_{n,k},$$
and taking into consideration that
$$L_n^\alpha(0)=\binom{n+\alpha}{n}=w(\alpha,n)=:w(n),$$
we denote by
$$R_n^\alpha(x)=\frac{L_n^\alpha(x)}{\binom{n+\alpha}{n}}.$$
Thus $R_n^\alpha(0)=1$ and
\begin{equation}\label{Rnorm} \left\|\sqrt{\frac{w(n)}{\Gamma(\alpha+1)}}R_n^\alpha\right\|_{2,\alpha}=1, \ws \ws \|\tilde{R}_n^\alpha\|_\infty=\tilde{R}_n^\alpha(0)=1.\end{equation}
With this notation we have that
\begin{equation}\label{rr}T_t^\alpha(R_n^\alpha,x)=R_n^\alpha(x)R_n^\alpha(t),\end{equation}
see e.g. \cite{gm}.

According to \cite[Theorem 1]{gm} for all $\alpha\ge 0$ and  $1\le p \le \infty$, considering the translated function as a function of the variable $x$,
\begin{equation}\label{trnorm}\|e^{-\frac{t}{2}}T_t^\alpha(f,x)\|_{p,\alpha}\le \|f\|_{p,\alpha}.\end{equation}

The corresponding convolution is
\begin{equation}\label{colag}(f*g)(t):=\int_0^\infty T_t^\alpha(f,x)g(x)e^{-x}x^\alpha dx.\end{equation}
Again by  \cite[Theorem 1]{gm}, for all $\alpha\ge 0$ and  $1\le p, q, r \le \infty$ with $\frac{1}{r}=\frac{1}{p}+\frac{1}{q}-1$; if $f\in L^p_\alpha$ and $g \in L^q_\alpha$, then
\begin{equation}\|f*g\|_{r,\alpha} \le \|f\|_{p,\alpha}\|g\|_{q,\alpha}.\end{equation}

The first theorem of this section is a Kolmogorov -Riesz type theorem, where the standard equicontinuity property is replaced by the one based on Laguerre translation.

\medskip

\begin{theorem}\label{lpr}  Let $1\le p <\infty$ and $\alpha \ge 0$. A bounded set $K\subset L^p_\alpha$ is precompact if and only if the properties below are fulfilled.\\
${\bf P_a}.$ For all $\varepsilon >0$ there is an $R>0$ such that  for all $f\in K$
\begin{equation}\label{farok} \left(\int_R^\infty |\tilde{f}|^pd\mu_\alpha\right)^{\frac{1}{p}}<\varepsilon. \end{equation}
${\bf P_b.}$ For all $\varepsilon$ and $M_0$ positive numbers there is a $\delta>0$ such that for all $0\le t \le M_0$, $0\le h \le \delta$ and $f\in K$
\begin{equation}\label{leqf}\left(\int_0^\infty \left|\left(T_{t+h}^\alpha f(x)-T_{t}^\alpha f(x)\right)e^{-\frac{x}{2}}\right|^pd\mu_\alpha(x)\right)^{\frac{1}{p}} <\varepsilon. \end{equation}

\end{theorem}

\proof

First we assume that $K$ is precompact. Then for an $\varepsilon>0$ there is a finite $\frac{\varepsilon}{2}$-net, $H_{\frac{\varepsilon}{2}}=\{u_1,\dots , u_j\}$ in $K$ $(j=j(\varepsilon))$. Since $C_0$ - the continuous functions with compact support - are dense in $L_{p,\alpha}$, there is an $S_\varepsilon \subset C_0$, such that $S_\varepsilon =\{\Phi_1, \dots ,\Phi_j\}$ and $\|u_l-\Phi_l\|_{p,(\alpha)}<\frac{\varepsilon}{2}$ for all $u_l \in H_{\frac{\varepsilon}{2}}$. That is there is an $R=R_\varepsilon >0$ such that the closed ball $\overline{B(0,R)}$ contains the support of each $\Phi \in S_\varepsilon$ wich ensures  ${\bf P_a}$.
To prove ${\bf P_b}$ take an $f\in L^p_\alpha$, $t \in [0,M_0]$ with some $M_0$, and let $\Phi$ be the element of $S_\varepsilon$ closest to $f$ with $\mathrm{supp} \Phi \in (0,R)$. In view of \eqref{transl}
$$\left(\int_0^\infty \left|\left(T_{t+h}^\alpha \Phi(x)-T_{t}^\alpha \Phi(x)\right)e^{-\frac{x}{2}}\right|^pd\mu_\alpha(x)\right)^{\frac{1}{p}}$$ $$=\left(\int_0^\infty \left|c_\alpha\int_0^\pi \left(\Phi(x+t+h+2\sqrt{x(t+h)}\cos\vartheta)e^{-\sqrt{x(t+h)}\cos\vartheta}j_{\alpha-\frac{1}{2}}(\sqrt{x(t+h)}\sin\vartheta)\right.\right.\right.$$ $$\left.\left.\left.-\Phi(x+t+2\sqrt{xt}\cos\vartheta)e^{-\sqrt{xt}\cos\vartheta}j_{\alpha-\frac{1}{2}}(\sqrt{xt}\sin\vartheta)\right)\sin^{2\alpha}\vartheta d\vartheta e^{-\frac{x}{2}}\right|^p x^\alpha dx\right)^{\frac{1}{p}}=(*).$$
For sake of simplicity let us denote by
$$\Psi(x,t,\vartheta):= \Phi(x+t+2\sqrt{xt}\cos\vartheta)e^{-\sqrt{xt}\cos\vartheta}j_{\alpha-\frac{1}{2}}(\sqrt{xt}\sin\vartheta).$$
Recalling that $\mathrm{supp}\Phi \subset (0,R)$ and $t\in [0, M_0]$, if $R<|\sqrt{x}-\sqrt{t}|$ then $\Psi(x,t,\vartheta)=0$, that is if $\sqrt{x}>\sqrt{M_0}+R$ then $\Psi(x,t,\vartheta)=0$. Considering that difference of the arguments of the three functions above are at most $2\sqrt{xh}\le c(R, M_0)\sqrt{h}$ and $\Psi(x,t,\vartheta)$ is compactly supported and continuous, there is a $\delta>0$ such that if $0\le h\le \delta$
$$\left|\Psi(x,t+h,\vartheta)-\Psi(x,t,\vartheta)\right|\le\varepsilon$$
and so
$$(*)\le \varepsilon \left(\int_0^\infty e^{-p\frac{x}{2}} x^\alpha dx\right)^{\frac{1}{p}}.$$
Then, because the set of $\Phi$-s in question is finite we can choose $R$ and $\delta$ uniformly. Since the norm of translation is bounded on $[0,M_0]$, cf. \eqref{trnorm} (and it is a linear operator), we can finish this part with triangle inequality.

Vice versa assuming ${\bf P_a}$ and ${\bf P_b}$ let
$$V_af(x):=\frac{1}{A}\int_0^aT_t^\alpha f(x)e^{-\frac{t}{2}}t^\alpha dt$$
if $0\le x \le R$ with some finite $R$, and we define $V_af(x)=0$ if $x>R$.
where $A=\int_0^a e^{-\frac{t}{2}}t^\alpha dt$. Then, applying H\"older's inequality and the symmetry of translation
\begin{equation}\label{maeq}|V_af(x+u)-V_af(x)|\le \frac{1}{A}\int_0^a|T_t^\alpha f(x+u)-T_t^\alpha f(x)|e^{-\frac{t}{2}}t^\alpha dt\end{equation} $$\le \frac{1}{A^{\frac{1}{p}}}e^{\frac{a}{2p'}}\left(\int_0^\infty|T_{x+u}^{\alpha}f(t)-T_x^{\alpha}f(t)|^pe^{-\frac{pt}{2}}t^\alpha dt\right)^{\frac{1}{p}},$$
and similarly
\begin{equation}\label{V}|V_af(x)|\le \frac{1}{A^{\frac{1}{p}}}e^{\frac{a}{2p'}}\left(\int_0^\infty|T_x^{\alpha}f(t)|^pe^{-\frac{pt}{2}}t^\alpha dt\right)^{\frac{1}{p}}.\end{equation}
Let
$$F_{a,R}:=\left\{V_af(x): f \in K, \ws x \le R\right\}.$$
By \eqref{V} and \eqref{trnorm}
$$|e^{-\frac{x}{2}}V_af(x)|\le \frac{1}{A^{\frac{1}{p}}}e^{\frac{a}{2p'}}\|f\|_{p,\alpha}.$$
Thus with fixed $a$ and $R$, $F_{a,R}$ is bounded (by $M \frac{1}{A^{\frac{1}{p}}}e^{\frac{a}{2p'}+\frac{R}{2}}$ if $K$ is bounded by $M$) and (choosing $M_0=R$) by \eqref{maeq} and assumption ${\bf P_b}$ it is equicontinuous. Thus for an arbitrary $\varepsilon>0$ there is an $\varepsilon$-net, $V_af_1, \dots, V_af_n$, in $F_{a,R}$ such that $f_i \in K$, $i=1, \dots , n$.\\
Let $f\in K$ be arbitrary, and for an $\varepsilon>0$ choose $R$ according to property ${\bf P_a}$. Then again by H\"older and Fubini theorems
$$\| V_af-f\|_{p,\alpha}$$ $$\le \left(\int_0^R\left|e^{-\frac{x}{2}}\int_0^a\frac{1}{A}\left(T_t^{\alpha}f(x)-f(x)\right)e^{-\frac{t}{2}}t^\alpha dt\right|^px^\alpha dx\right)^{\frac{1}{p}}+\left(\int_R^\infty\left|e^{-\frac{x}{2}}f(x)\right|^px^\alpha dx\right)^{\frac{1}{p}}$$
$$\le \varepsilon +e^{\frac{a}{2p'}}\left(\frac{1}{A}\int_0^\infty e^{-\frac{px}{2}}\int_0^a\left|T_t^{\alpha}f(x)-f(x)\right|^pe^{-\frac{t}{2}}t^\alpha dt x^\alpha dx\right)^{\frac{1}{p}}$$ $$\le \varepsilon +e^{\frac{a}{2p'}}\left(\frac{1}{A}\int_0^ae^{-\frac{t}{2}}t^\alpha \int_0^\infty e^{-\frac{px}{2}}\left|T_t^{\alpha}f(x)-f(x)\right|^px^\alpha dx dt\right)^{\frac{1}{p}}$$ $$\le \varepsilon +e^{\frac{a}{2p'}}\sup_{0\le t \le a}\left(\int_0^\infty e^{-\frac{px}{2}}\left|T_t^{\alpha}f(x)-f(x)\right|^px^\alpha dx\right)^{\frac{1}{p}}.$$
Thus choosing $a$ according to property ${\bf P_b}$ small enough (assume that $a<1$), $\| V_af-f\|_{p,\alpha}<2\varepsilon$. \\
With these chosen $a$ and $R$ construct $F_{a,R}$ and for the previous $\varepsilon$ and $f$, from the $\varepsilon$-net in question choose $V_af_i$ such that $|V_af(x)-V_af_i(x)|<\varepsilon$ on $[0, R]$. Then we have
$$\| V_af-V_af_i\|_{p,\alpha}= \left(\int_0^R\left|e^{-\frac{x}{2}}( V_af(x)-V_af_i(x))\right|^px^\alpha dx\right)^{\frac{1}{p}}
$$ $$\le \varepsilon \left(\int_0^\infty e^{-\frac{px}{2}}x^\alpha dx\right)^{\frac{1}{p}}\le\varepsilon c(\alpha,p).$$
That is by the triangle inequality $\{f_i\}_{i=1}^n$ is a $(4+c(\alpha,p))\varepsilon$-net in $K$.

\medskip

\subsection{Laguerre translation on sequences}

Laguerre translation on sequences, to the best of my knowledge, have not appeared explicitly just as convolution of sequences, see e.g. \cite{ag}. The corresponding algebras are investigated in \cite{ka}. Below we derive the translation from convolution and investigate its properties.

Let $\alpha>-1$. Recalling the notation above we introduce the discrete weights
$$w(k)=w(\alpha,k)=\binom{k+\alpha}{k},$$
and the space of real sequences, $a=\{a(k)\}_{k=0}^\infty$. For $1\le p<\infty$ $a\in l^p_\alpha$ if
$$\|a\|_{p,\alpha} :=\left(\sum_{k=0}^\infty|a(k)|^pw(k)\right)^{\frac{1}{p}}<\infty.$$
$$\|a\|_{\infty,\alpha}=\|a\|_{\infty}$$
independently of $\alpha$.

If $a\in l^p_\alpha$ and $b\in l^{p'}_\alpha$ we denote by
$$\langle a,b\rangle =\sum_{k=0}^\infty a(k)b(k)w(k).$$

\remark

Certainly, in discrete case the criterion of precompactness is simpler, cf. \cite[Theorem 4]{hoh}. It is as follows.

For any $1\le p <\infty$ and $\alpha \ge 0$ a set $K \subset l^p_\alpha$ is precompact if and only if it is pointwise bounded (i.e. for all $n\in\mathbb{N}$ there is an $M(n)$ such that for all $a\in K $ $a(n)\le M(n)$), and the next property fulfils.

\medskip

${\bf P_{as}.}$ For all $\varepsilon >0$ there is an $N\in\mathbb{N}$ such that for all $a\in K $
 \begin{equation}\label{afarok}\left(\sum_{k=N+1}^\infty |a(k)|^pw(k)\right)^{\frac{1}{p}} <\varepsilon. \end{equation}

Indeed, suppose, that $K$ is precompact. Then it is obviously bounded in $l^p_\alpha$. Since $w(k)\ge 1$, it is  pointwise bounded as well. For $\frac{\varepsilon}{2}$ it has a finite $\frac{\varepsilon}{2}$-net, $b_1, \dots, b_n$, say. Since the finite sequences are dense in $l^p_\alpha$, there are $s_i$ finite sequences, such that $\|p_i-s_i\|_{p,\alpha}<\frac{\varepsilon}{2}$, $i=1, \dots, n$. Thus the maximal length of $s_i$-s is an appropriate choice of $N$.

Assume now, that $K$ is pointwise bounded and fulfils ${\bf P_{as}}$. Choose an $N$  for an arbitrary $\varepsilon$ ensured by ${\bf P_{as}}$. Let $S_N:=\{a^N:=a(0), \dots, a(N) : a\in K\}$ be the set of the $N+1$-long initial parts of the sequences in $K$. Then the distance of $K$ and $S_N$ is at most $\varepsilon$. Because $S_N$ is finite dimensional it is bounded in $l^p_\alpha$ ($S_N$ is bounded by $cM_N N^{\alpha+1}$, where $M_N:=\max_{n=1}^N M(n)$) and (again by being finite dimensional) $S_N$ contains a finite $\varepsilon$-net, and the corresponding sequences form a $2\varepsilon$-net in $K$.

\medskip

To define translation on the spaces of sequences, our starting point is the convolution defined in \cite{ag}.\\ Let $\alpha >\alpha_0=\frac{-5+\sqrt{17}}{2}$. Then by \cite[Theorem 1]{ag}
\begin{equation}\label{gam}\gamma(n,m,k):=\gamma(n,m,k,\alpha)=\int_0^\infty R_n^\alpha(x)R_m^\alpha(x)R_k^\alpha(x)e^{-2x}x^\alpha dx >0,\end{equation}
for all $k, m, n \in \mathbb{N}$. In view of \cite[(4.2)]{ag}
\begin{equation}\label{sgam}\sum_{k=0}^\infty \gamma(n,m,k)w(k)=1\end{equation}
for all $m, n \in \mathbb{N}$.

Let $a \in l^p_\alpha$, $b \in l^q_\alpha$, where $\frac{1}{p}+\frac{1}{q}\ge 1$. In \cite[(4.5)]{ag} the next convolution is defined.
\begin{equation}\label{aconv}(a*b)(k):=\sum_{m=0}^\infty \sum_{n=0}^\infty a(m)b(n)\gamma(n,m,k)w(n)w(m).\end{equation}
Similarly to the $L^p$-cases if $1\le p,q,r \le \infty$ and $\frac{1}{r}=\frac{1}{p}+\frac{1}{q}-1$,
\begin{equation}\label{cn}\|a*b\|_{r,\alpha} \le \|a\|_{p,\alpha}\|b\|_{q,\alpha},\end{equation}
see \cite[(4.6)]{ag}.

Thus, fixing $\alpha>\alpha_0$, we can define the corresponding translation as
\begin{equation}\label{selt} T_k(a)(n):=\sum_{m=0}^\infty  a(m)w(m)\gamma(n,m,k).\end{equation}
By symmetry we immediately have that
\begin{equation} T_k(a)(n)=T_n(a)(k).\end{equation}
So \eqref{aconv} can be written as
\begin{equation}(a*b)(k)=\sum_{n=0}^\infty b(n)T_k(a)(n)w(n)=\langle b, T_k(a)\rangle =\langle a, T_k(b)\rangle.\end{equation}

As on functions, Laguerre translation on sequences is also a bounded operator.

\begin{proposition}\label{Tsn} Let $\alpha>\alpha_0$ and $1\le p \le \infty$ Then for all $k$
$$\|T_k(a)\|_{p,\alpha}\le \|a\|_{p,\alpha}.$$
\end{proposition}

\proof
Let $a \in l^\infty$. In view of \eqref{selt}, \eqref{gam} and \eqref{sgam}
$$|T_k(a)(n)|=|\sum_{m=0}^\infty  a(m)w(m)\gamma(n,m,k)|\le \|a\|_\infty \sum_{m=0}^\infty \gamma(n,m,k)w(m)=\|a\|_\infty.$$
If $a \in l^1_\alpha$
$$\sum_{n=0}^\infty |T_k(a)(n)|w(m)=\sup_{b, \|b\|_\infty \le 1}\sum_{n=0}^\infty b(n)T_k(a)(n)w(n)$$ $$=\sup_{b, \|b\|_\infty \le 1}\sum_{m=0}^\infty \sum_{n=0}^\infty a(m)b(n)\gamma(n,m,k)w(n)w(m)=\sup_{b, \|b\|_\infty \le 1}(a*b)(k).$$ According to \eqref{cn}
$$ \|a*b\|_\infty \le  \|a\|_{1,\alpha},$$
which implies the statement in $l^1_\alpha$. Finally interpolation ensures the result.

\medskip

 To prepare the next subsection, we mention the standard connection between the corresponding spaces of functions and sequences.\\
 For a function $f$ in some $L^p_\alpha$ let us denote the corresponding sequence by $\hat{f}=a_f$, where
\begin{equation} \hat{f}(n):=\int_0^\infty \tilde{f}\tilde{R}_n^\alpha d\mu_\alpha, \end{equation}
 for $n \in \mathbb{N}$.\\
 Let $a$ be a sequence. The corresponding function can be defined as $\check{a}=f_a(x)$, where
 $$\tilde{f}_a \sim \sum_{m=0}^\infty  a(m)w(m)\tilde{R}_m^\alpha$$
 if the series is convergent in some sense.\\
 Let $1\le p \le 2$. Then, as usual, if $f\in L^p_\alpha$ then $\hat{f} \in l^{p'}_{\alpha}$, and if $a\in l^{p}_{\alpha}$ then $\check{a}\in L^{p'}_{\alpha}$ and the operators map $L^p_\alpha$ to $l^{p'}_{\alpha}$ and $l^p_\alpha$ to $L^{p'}_{\alpha}$ are bounded.\\
 Indeed, by interpolation it is a consequence of Parseval's formula and the next inequalities:
 $$|\hat{f}(k)|=\left|\int_0^\infty \tilde{f}\tilde{R}^{\alpha}_kd\mu_\alpha\right|\le \|f\|_{1,\alpha},$$
and
$$\|\tilde{f}_a\|_\infty \le\sum_{m=0}^\infty  |a(m)|w(m)\|\tilde{R}_m^{\alpha}\|_\infty =\|a\|_{1,\alpha},$$
where \eqref{Rnorm} is considered.

 Let $a$ and $b$ are in $l^1_\alpha$, say, then the series $\tilde{f}_a$ and $\tilde{g}_b$  are uniformly convergent.  Thus by \eqref{aconv}
 \begin{equation}\tilde{f}_a(x)\tilde{g}_b(x)e^{-\frac{x}{2}}=\sum_{k=0}^\infty  (a*b)(k)w(k)\tilde{R}_k^\alpha,\end{equation}
 that is
 \begin{equation}\label{sz}\check{a}\check{b}\sim a*b,\end{equation}
 see \cite{ag}.

 Similarly \eqref{aconv} ensures that if $\hat{f}=a_f$ then
 \begin{equation}\label{Tnk}\widehat{\tilde{f}\tilde{R}_n^{\alpha}}(k)=T_n(a_f)(k).\end{equation}

 Moreover  \eqref{rr} implies that
 \begin{equation}\label{38}\widehat{f*g}(n)=\hat{f}(n)\hat{g}(n), \ws \ws   \widehat{ T_t^\alpha(f)}(n)=\hat{f}(n)R_n^\alpha(t),\end{equation}
 provided that $\alpha \ge 0$, $f\in L^p_\alpha$ and $g \in L^q_\alpha$ and  $\frac{1}{p}+\frac{1}{q}\ge 1$.

\subsection{Pego-type theorem with Laguerre transformation} At first we introduce the notion of equicontinuity in mean with respect to sequences.

\medskip

\noindent ${\bf P_{bs}}.$ A set $K \subset  l^p_\alpha$ is equicontinuous in mean if for all $\varepsilon>0$ there is an $N \in \mathbb{N}$ such that for all $j>N$ and $a \in K$
\begin{equation}\label{aeqc}  \left(\sum_{k=0}^\infty |T_j(a)(k)-a(k)|^pw(k)\right)^{\frac{1}{p}}<\varepsilon.\end{equation}

\medskip

Subsequently we need the the next extra property with respect to functions.

\medskip

\noindent ${\bf P_{a0}}.$ A set $K \subset  L^p_\alpha$ is equivanishing at zero if for all $\varepsilon>0$ there is a $\delta>0$ such that for all $f\in K$
\begin{equation}\label{eqv0}\left(\int_0^\delta |\tilde{f}|^pd\mu_\alpha\right)^{\frac{1}{p}}\le \varepsilon.\end{equation}

\medskip

After this preparation we are in position to state a Pego-type theorem.

\medskip

\begin{theorem}\label{lL}  Let $1\le p \le 2$ and $\alpha \ge 0$. \\
$(\mathrm{a})$ If $K \subset  l^p_\alpha$ is bounded and fulfils ${\bf P_{as}}$, then $\check{K} \subset L^{p'}_\alpha$ fulfils ${\bf P_{b}}$.\\
$(\mathrm{b})$ If $K \subset  L^p_\alpha$ fulfils ${\bf P_{b}}$, then $\hat{K} \subset l^{p'}_\alpha$ fulfils ${\bf P_{as}}$.\\
$(\mathrm{c})$ If $K \subset  l^p_\alpha$ fulfils ${\bf P_{bs}}$, then $\check{K} \subset L^{p'}_\alpha$ fulfils ${\bf P_{a}}$.\\
$(\mathrm{d})$ If in addition $\alpha >\frac{1}{2}$, $K \subset  L^p_\alpha$ is bounded and fulfils ${\bf P_{a}}$ and ${\bf P_{a0}}$, then $\hat{K} \subset l^{p'}_\alpha$ fulfils ${\bf P_{bs}}$.
\end{theorem}

\proof
$(\mathrm{a})$: Let $f=f_a \in \check{K} \subset L^{p'}_\alpha$. Then by \eqref{38} and \eqref{Rnorm}
$$\left(\int_0^\infty\left|\left(T_{t+h}^{\alpha}f(x)-T_t^{\alpha}f(x)\right|e^{-\frac{x}{2}}\right)^{p'}d\mu_\alpha(x)\right)^{\frac{1}{p'}}$$ $$\le c\left(\sum_{k=0}^\infty\left|a(k)\left(R_k^{\alpha}(t+h)-R_k^{\alpha}(t)\right)\right|^pw(k)\right)^{\frac{1}{p}}$$ $$\le c\left(\sum_{k=0}^N\left|a(k)\left(R_k^{\alpha}(t+h)-R_k^{\alpha}(t)\right)\right|^pw(k)\right)^{\frac{1}{p}}+2ce^{\frac{t}{2}}\left(\sum_{k=N+1}^\infty |a(k)|^pw(k)\right)^{\frac{1}{p}}$$ $$=S_1+S_2.$$
Since $\left\|\left(R_k^{\alpha}\right)'\right\|_\infty \le ck \left\|R_{k-1}^{(\alpha+1)}\right\|_\infty $, see \cite[(5.1.14)]{sz},
$$S_1\le cNMe^{\frac{t}{2}}h,$$
where $K$ is bounded by $M$. Thus, recalling that $t\in [0,M_0]$ by ${\bf P_{as}}$ if $N$ is large enough,
$S_2\le \frac{\varepsilon}{2}$, and we can choose $\delta$ so small such that $S_1\le \frac{\varepsilon}{2},$ if $0\le h \le \delta$.\\
$(\mathrm{b})$: Let $f\in K$. Take a function $g\in C_0$ such that\\ $g\ge 0$, $\mathrm{supp}g\subset [0,\delta]$ and $\int_0^\infty g(x)e^{-x}d\mu_\alpha(x)=1$. As $\hat{g}(k) \to 0$, we can choose $N$ so large that $| \hat{g}(k)|<\frac{1}{2}$ if $k\ge N$. Thus applying \eqref{38}, \eqref{colag}, Minkowski inequality and the symmetry of translation we have
$$\left(\sum_{k=N}^\infty |\hat{f}(k)|^{p'}w(k)\right)^{\frac{1}{p'}}\le 2 \left(\sum_{k=N}^\infty |\hat{f}(k)(1-\hat{g}(k))|^{p'}w(k)\right)^{\frac{1}{p'}}$$ $$\le c\left(\int_0^\infty\left|(f(t)-f*g(t))e^{-\frac{t}{2}}\right|^pd\mu_\alpha(t)\right)^{\frac{1}{p}}$$ $$=c\left(\int_0^\infty\left|\left(f(t)-\int_0^\infty T_t^{\alpha}f(x)g(x)e^{-x}d\mu_\alpha(x)\right)e^{-\frac{t}{2}}\right|^pd\mu_\alpha(t)\right)^{\frac{1}{p}}$$ $$\le c\int_0^\infty\left(\int_0^\infty \left(|f(t)-T_x^{\alpha}f(t)|e^{-\frac{t}{2}}\right)^pd\mu_\alpha(t)\right)^{\frac{1}{p}}g(x)e^{-x}d\mu_\alpha(x)$$ $$\le c\sup_{0\le x \le \delta}\|f-T_x^{\alpha}f\|_{p,\alpha}.$$
$(\mathrm{c})$: Let $N$ be the appropriate index for $\varepsilon$ in ${\bf P_{bs}}$, and let us define a finite sequence $b=(b(N), \dots, b(N+n))$ such that $b_k\ge 0$, $k=N,\dots, N+n$ and $\sum_N^{N+n}b(k)w(k)=1$. Let us denote by  $\check{b}=g$ and for an $a\in K$ $\check{a}=f$. For $\varepsilon>0$ we can take an apropriate initial part of $a$, $a_L:=(a(0),a(1), \dots, a(L))$, and $\check{a}_L=p_L$ such that $\|f-p_L\|_{p',\alpha}\le \varepsilon $.\\
Since $\lim_{x\to\infty}g(x)e^{-x}=0$, there is an $R$ such that $|g(x)e^{-x}|\le \frac{1}{2}$ if $x>R$. Considering that $|g(x)e^{-x}|\le 1$ on $x\ge 0$, in view of \eqref{sz}, as above
$$\left(\int_R^\infty|\tilde{f}(x)|^{p'}d\mu_\alpha(x)\right)^{\frac{1}{p'}} \le \varepsilon+\left(\int_R^\infty|\tilde{p_L}(x)|^{p'}d\mu_\alpha(x)\right)^{\frac{1}{p'}}$$ $$\le \varepsilon+ 2\left(\int_R^\infty|\tilde{p_L}(x)(1-g(x)e^{-x})|^{p'}d\mu_\alpha(x)\right)^{\frac{1}{p'}}$$ $$=\varepsilon+\left(\int_R^\infty|\tilde{p_L}(x)-\tilde{p_L}(x)g(x)e^{-x}|^{p'}d\mu_\alpha(x)\right)^{\frac{1}{p'}}\le \varepsilon+c\|a_L-a_L*b\|_{p,\alpha}$$ $$=\varepsilon+c\left(\sum_{k=0}^L\left|a(k)-\sum_{j=N}^{N+n}b(j)T_k(a)(j)w(j)\right|^pw(k)\right)^{\frac{1}{p}}$$ $$\le\varepsilon+c\left(\sum_{k=0}^\infty\left|\sum_{j=N}^{N+n}\left[a(k)-T_j(a)(k)\right]b(j)w(j)\right|^pw(k)\right)^{\frac{1}{p}}$$ $$\varepsilon+\le c \sup_{j\ge N}\left(\sum_{k=0}^\infty |a(k)-T_j(a)(k)|^pw(k)\right)^{\frac{1}{p}}\le(1+c)\varepsilon.$$
$(\mathrm{d})$: According to \eqref{Tnk} and \eqref{Rnorm}
$$\left(\sum_{k=0}^\infty |a(k)-T_j(a)(k)|^{p'}w(k)\right)^{\frac{1}{p'}}=\left(\sum_{k=0}^\infty \left|\hat{f}(k)-\widehat{\tilde{R}_j^{\alpha}\tilde{f}}(k)\right|^{p'}w(k)\right)^{\frac{1}{p'}}$$ $$\le c\left(\int_0^\infty \left|f(x)\left(1- R_j^{\alpha}(x)e^{-x}\right)e^{-\frac{x}{2}}\right|^pd\mu_\alpha(x))\right)^{\frac{1}{p}}\le 2c\left(\int_0^\delta \left|f(x)e^{-\frac{x}{2}}\right|^pd\mu_\alpha(x))\right)^{\frac{1}{p}}$$ $$+ c\left(\int_\delta^R \left|f(x)\left(1- R_j^{\alpha}(x)e^{-x}\right)e^{-\frac{x}{2}}\right|^pd\mu_\alpha(x))\right)^{\frac{1}{p}}$$ $$+2c\left(\int_R^\infty \left|f(x)e^{-\frac{x}{2}}\right|^pd\mu_\alpha(x))\right)^{\frac{1}{p}}=I+II+III.$$
By the assumption $I+III\le 4c\varepsilon$. To estimate $II$ we consider that
$$\left(R_j(x)e^{-x}\right)'=e^{-x}\frac{-L_{j-1}^{(\alpha+1)}(x)-L_j^{\alpha}(x)}{\binom{j+\alpha}{j}}=\frac{-L_{j}^{(\alpha+1)}(x)e^{-x}}{\binom{j+\alpha}{j}},$$ see \cite[(5.1.13)]{sz}. In view of \cite[(2.8)]{m}
\begin{equation}\label{infn}\tilde{R}_n^{\alpha}(x)x^{\frac{\alpha}{2}}\le \frac{C}{n^{\frac{\alpha}{2}}} \left\{\begin{array}{ll}(x\nu){\frac{\alpha}{2}}, \ws \ws 0\le x\le \frac{1}{\nu}\\
(x\nu)^{-\frac{1}{4}}, \ws \ws \frac{1}{\nu}< x \le \frac{\nu}{2}\\
\left(\nu\left(\nu^{\frac{1}{3}}+|x-\nu|\right)\right)^{-\frac{1}{4}}, \ws \ws \frac{\nu}{2}< x \le \frac{3\nu}{2}\\
e^{-\gamma x}, \ws \ws \frac{3\nu}{2}<x, \end{array}\right.\end{equation}
where $\nu=4n+2\alpha+2$. Thus
\begin{equation}\label{derb}\frac{\left|L_{j}^{(\alpha+1)}(x)\right|e^{-x}}{\binom{j+\alpha}{j}}\le c\frac{j^{\frac{1}{4}-\frac{\alpha}{2}}}{\delta^{\frac{\alpha+1}{2}+\frac{1}{4}}}, \ws\ws x\in(\delta,R).\end{equation}
Thus, denoting the bound of $K$ with $M$
$$II\le c(\delta)Rj^{\frac{1}{4}-\frac{\alpha}{2}}\|f\|_{p,\alpha}\le c(\delta,R,M)j^{\frac{1}{4}-\frac{\alpha}{2}},$$
which is small if $j$ is large enough.

\medskip

The computation in proof of (d) ensures the regular behavior of translation on sequences.

\medskip

\begin{proposition} Let $\alpha >\frac{1}{2}$, $1\le p <\infty$. If $K\subset l^p_\alpha$ is precompact, then $K$ fulfils ${\bf P_{bs}}$.\end{proposition}

\proof
By property ${\bf P_{as}}$ it is enough to prove the statement for finite sets of finite sequences.
Let $S:=s_1, \dots, s_n$ and  $s_i=\{a_{ik}\}_{k=0}^{n_i}$ and $p_i=\sum_{k=0}^{n_i}a_{ik}R_k^\alpha$. For any $\varepsilon>0$ there are $\delta_i$ and $R_i$ such that $\left(\int_0^{\delta_i}|\tilde{p}_i|^2d\mu_\alpha\right)^{\frac{1}{2}}\le\varepsilon$ and $\left(\int_{R_i}^{\infty}|\tilde{p}_i|^2d\mu_\alpha\right)^{\frac{1}{2}}\le\varepsilon$, $i=1, \dots, n$. Let $R:=R(\varepsilon_1,\varepsilon)=\max_{i=1}^n R_i$, $N:=N(\varepsilon_1)=\max_{i=1}^n n_i$, $\delta:=\delta(\varepsilon_1,\varepsilon)=\min_{i=1}^n \delta_i$, $M:=M(\varepsilon_1)=\max_{i=1}^n\|p_i\|_{2,\alpha}$. With the abbreviation $a_{ik}=\hat{p}_i(k)$, as above we have
$$\left|\hat{p}_i(k)-T_j(\hat{p}_i)(k)\right|=\left|\int_0^\infty \left(p_i(x)-R_j^{\alpha}(x)p_i(x)e^{-x}\right)R_k(x)e^{-x}x^\alpha dx\right|$$ $$\le \|R_k^{\alpha}\|_{2,\alpha}\left(2\left(\int_0^\delta |\tilde{p}_i|^2 d\mu_\alpha \right)^{\frac{1}{2}}+\|p_i\|_{2,\alpha}\sup_{[\delta,R]}\left|1-R_j(x)e^{-x}\right|+2\left(\int_R^\infty |\tilde{p}_i|^2 d\mu_\alpha \right)^{\frac{1}{2}}\right).$$
By \eqref{derb} and \eqref{Rnorm}
$$\left|\hat{p}_i(k)-T_j(\hat{p}_i)(k)\right|\le ck^{-\frac{\alpha}{2}} \left(4\varepsilon+ M R \frac{j^{\frac{1}{4}-\frac{\alpha}{2}}}{\delta^{\frac{\alpha+1}{2}+\frac{1}{4}}}\right)\leq 5c\varepsilon k^{-\frac{\alpha}{2}}$$
if $j$ is large enough. Thus
$$\|\hat{p}_i-T_j(\hat{p}_i)\|_{p,\alpha}\le 5c\varepsilon \left(\sum_{k=0}^Nk^{-p\frac{\alpha}{2}}w(k)\right)^{\frac{1}{p}}\le 5c\varepsilon N^{\frac{\alpha+1}{p}-\frac{\alpha}{2}}.$$
Since $N$ is independent of $\varepsilon$, it is arbitrary small if $j$ is large enough. According to Proposition \ref{Tsn} the proof can be finished with a triangle inequality.

\medskip

The main theorem of this section os the next corollary.

\medskip

\begin{cor} Let $\alpha \ge 0$. A set $K\subset L^2_\alpha$ is precompact if and only if it satisfies ${\bf P_{b}}$, and the corresponding set $\hat{K}$ is pointwise bounded.
\end{cor}

\proof
If $K$ is precompact in $L^2_\alpha$, Theorem \ref{lpr} ensures equicontinuity in mean. Moreover $K$ is bounded in  $L^2_\alpha$, so $\hat{K}$ is bonded in $l^2_\alpha$ and as above, pointwise too.\\
If $K$  is equicontinuous in mean, by Theorem \ref{lL} $\hat{K}$ is equivanishing in $l^2_\alpha$, and the pointwise boundedness ensures precompactness of $\hat{K}$ and of $K$ too.

\medskip

\section{Bessel translation method}

Bessel translation and Hankel transformation are widely examined by several authors. We mention here only a few examples. One of our main sources is an early paper, \cite{le}. By Bessel translation modulus of smoothness and the related best approximation can be investigated, see \cite{p}, and Nikol'skii inequalities for entire functions can be proved, see \cite{abdh}. For Fourier-Bessel transformation, similarly to the standard Fourier transformation, uncertainty results  are derived, see \cite{gj}, which are somehow in concordance with Pego-type theorems. Hereinafter we concentrate to compactness criteria.

Let $L_{p,\alpha}$ be the space of measurable functions on $\mathbb{R}_+$ equipped with the norm
$$\|f\|_{p,(\alpha)}=\left(\int_0^\infty|f(x)|^p x^{2\alpha+1}dx\right)^{\frac{1}{p}}, \ws \ws \alpha \ge -\frac{1}{2}, \ws 1\le p <\infty, \ws \ws  \ws \|f\|_{\infty,(\alpha)}=\|f\|_{\infty, \mathbb{R}_+}.$$
The dual space of $L_{p,\alpha}$ is denoted by $L_{p',\alpha}$, where
$$\frac{1}{p}+\frac{1}{p'}=1.$$
First we introduce the Bessel translation of an integrable function, see \cite[(5.19)]{le}:
\begin{equation}\label{t} T_{t,\alpha}f(s)=\int_0^\pi f(\sqrt{t^2+s^2-2st\cos\varphi})d\mu(\varphi),\end{equation}
where $d\mu(\varphi)=c_\alpha \sin^{2\alpha}\varphi d\varphi$, and $c_\alpha=\left(\int_0^\pi \sin^{2\alpha}\varphi d\varphi\right)^{-1}$.
The symmetry of the definition implies
\begin{equation}\label{st} T_{t,\alpha}f(s)= T_{s,\alpha}f(t).\end{equation}

We need the entire Bessel functions which are
\begin{equation}\label{B}j_\alpha(z)=\Gamma(\alpha+1)\left(\frac{2}{z}\right)^\alpha J_\alpha(z)=\sum_{k=0}^\infty\frac{(-1)^k\Gamma(\alpha+1)}{\Gamma(k+1)\Gamma(k+\alpha+1)}\left(\frac{z}{2}\right)^{2k}.\end{equation}
Subsequently the next properties will be necessary. The norm of an entire Bessel function on the half-line is attained at zero:
\begin{equation}\label{j1}\|j_\alpha\|_{\infty, \mathbb{R}_+}=j_\alpha(0)=1.\end{equation}
The derivative of $j_\alpha$ can be expressed as
\begin{equation}\label{j2}j_\alpha'(z)=-\frac{1}{2(\alpha+1)}zj_{\alpha+1}(z),\end{equation}
see \cite{be}.

As it is mentioned in the introduction, similarly to the exponential case, the Bessel translation fulfils
$$T_{t,\alpha}j_\alpha(\lambda x)=j_\alpha(\lambda t)j_\alpha(\lambda x),$$
see \cite{le}.

The operator norm of Bessel translation is one, see \cite[(2.24)]{p}
\begin{equation}\label{opntr}\|T_tf\|_{p,(\alpha)}\le \|f\|_{p,(\alpha)}, \ws \ws 1\le p \le \infty.\end{equation}

The Kolmogorov-Riesz type theorem of this section is as follows.

\medskip

\begin{theorem}\label{best} $1\le p < \infty$ $K\subset L_{p,\alpha}$ is bounded. $K$ is precompact if and only if the two conditions below are fulfiled.\\
${\bf P_A.}$ $K$ is equivanishing, i.e.
\begin{equation}\label{c2} \forall \varepsilon>0 \ws \exists R>0, \ws \forall f\in K \ws \left(\int_R^\infty|f(x)|^px^{2\alpha+1}dx\right)^{\frac{1}{p}}<\varepsilon.\end{equation}
${\bf P_B.}$ $K$ is $L_{p,\alpha}$-equicontinuous, i.e. $\forall \varepsilon>0$ and $\forall$  $M_0>0$ there is a $\delta>0$ such that $\forall$ $0\le h \le \delta$, $\forall$ $f\in K$, $\forall$ $t\in[0,M_0]$
\begin{equation}\label{c3} \left(\int_0^\infty|T_{t+h,\alpha}f(s)-T_{t,\alpha}f(s)|^p s^{2\alpha+1}ds\right)^{\frac{1}{p}}<\varepsilon.\end{equation}
\end{theorem}

\proof
To prove property ${\bf P_A}$ we can proceed just as in the proof of Theorem \ref{lpr}.\\
To prove ${\bf P_B}$, with the notation of Theorem \ref{lpr} we take an $f\in L_{p,\alpha}$ and a $\Phi$ in $S_{\frac{\varepsilon}{3}}$ closest to $f$.\\
Denoting by $(s,t,\varphi):=(\sqrt{s^2+t^2-2st\cos\varphi})$, in view of \eqref{t}
$$\|T_{t+h,\alpha} \Phi- T_{t,\alpha}\Phi\|_{p,(\alpha)}$$ $$=\left(\int_0^\infty \left|\int_0^\pi \left(\Phi(s,t+h,\varphi)-\Phi(s,t,\varphi)\right)d\mu(\varphi)\right|^ps^{2\alpha+1}ds\right)^{\frac{1}{p}}$$ $$\le \int_0^\pi \left(\int_0^\infty \left|\Phi(s,t+h,\varphi)-\Phi(s,t,\varphi)\right|^ps^{2\alpha+1}ds\right)^{\frac{1}{p}}d\mu(\varphi).$$
Since $||s-t|-h|\le (s,t,\varphi),(s,t+h,\varphi)\le s+t+h$, recalling that $\Phi \in S_{\frac{\varepsilon}{3}}$ and $t<M_0$
$$\|T_{t+h,\alpha} \Phi- T_{t,\alpha}\Phi\|_{p,(\alpha)}\le \int_0^\pi \left(\int_0^R \left|\Phi(s,t+h,\varphi)-\Phi(s,t,\varphi)\right|^ps^{2\alpha+1}ds\right)^{\frac{1}{p}}d\mu(\varphi),$$
where $R=R_{\frac{\varepsilon}{3}}+M_0+1$.\\ 
$|\sqrt{s^2+(t+h)^2-2s(t+h)\cos\varphi}-\sqrt{s^2+t^2-2st\cos\varphi}|\le h$ (cf. \cite[(2.27)]{p}). Recall that there are finitely many uniformly continuous $\Phi$ in $S_{\frac{\varepsilon}{3}}$. Since they are equicontinuous (in standard sense) we can choose an appropriate $\delta$ to $\frac{\varepsilon}{3B}$, where $B=B(\varepsilon)= \left(\int_0^R s^{2\alpha+1}ds \right)^{\frac{1}{p}}$ such that $|\Phi(s,t+h,\varphi)-\Phi(s,t,\varphi)|<\frac{\varepsilon}{3B}$, thus considering that $\mu$ is a probability measure, $\|T_{t+h,\alpha} \Phi- T_{t,\alpha}\Phi\|_{p,(\alpha)}\le \frac{\varepsilon}{3}$. Considering the norm of the linear operator, \eqref{opntr}, by a triangle inequality \eqref{c3} is proved.

On the other hand supposing ${\bf P_A}$ and ${\bf P_B}$, we show that $K$ is precompact. To this we define
\begin{equation}M_{a}f(s):=\frac{1}{A}\int_0^a T_{t,\alpha}f(s)t^{2\alpha+1}dt,\end{equation}
where $A=\int_0^a t^{2\alpha+1}dt$.
By H\"older's inequality, then applying the symmetry of translation, cf. \eqref{st}
\begin{equation}\label{mfmf}|M_af(s+u)-M_af(s)|\le \frac{1}{A}\int_0^a|T_{t,\alpha}f(s+u)-T_{t,\alpha}f(s)|t^{2\alpha+1}dt\end{equation} $$\le \frac{1}{A^{\frac{1}{p}}}\left(\int_0^\infty |T_{t,\alpha}f(s+u)-T_{t,\alpha}f(s)|^pt^{2\alpha+1}dt\right)^{\frac{1}{p}}$$ $$=\frac{1}{A^{\frac{1}{p}}}\left(\int_0^\infty |T_{s+u,\alpha}f(t)-T_{s,\alpha}f(t)|^pt^{2\alpha+1}dt\right)^{\frac{1}{p}}.$$
Similarly, and recalling \eqref{opntr}
\begin{equation}\label{mfk}|M_af(s)|\le \frac{1}{A^{\frac{1}{p}}}\left(\int_0^\infty |T_{s,\alpha}f(t)|^pt^{2\alpha+1}dt\right)^{\frac{1}{p}}\le \frac{1}{A^{\frac{1}{p}}}\|f\|_{p,(\alpha)}.\end{equation}

Let $a$ and $R$ be fixed positive numbers, and let $F_{a,R}:= \{M_af(s): f\in K, s<R\}$. Then by \eqref{mfmf} and \eqref{c3} $F_{a,R}$ is equicontinuous (in standard sense), and by \eqref{mfk} and boundedness of $K$  $F_{a,R}$ is uniformly bounded. Thus by the theorem of Arzel\'a and Ascoli for all $\varepsilon>0$ there is a finite $\varepsilon$-net $N_\varepsilon \subset F_{a,R}$. \\
Let us denote the elements of $N_\varepsilon=\{M_af_1,\dots ,M_af_j\}$, where $f_i\in K$, $i=1,\dots , j$ and $j=j(\varepsilon)$, that is for all $f\in K$ there is an $f_i$, such that $|M_af(s)-M_af_i(s)|\le\varepsilon$ if $0\le s\le R$. We show that $\{f_i\}_{i=1}^j$ is an $\varepsilon$-net in $K$ if $\{M_af_i\}_{i=1}^j=N_{\frac{\varepsilon}{L}}$ is an $\frac{\varepsilon}{L}$-net in $F_{a,R}$ with a suitable $a$, $L$ and $R$ (will be given later).

Applying again H\"older's inequality and then Fubini's theorem
$$\|M_af-f\|_{p,(\alpha)}^p \le \int_0^\infty\left(\int_0^a\frac{1}{A}|T_{t,\alpha}f(s)-f(s)|t^{2\alpha+1}dt\right)^ps^{2\alpha+1}ds$$
$$\le \frac{1}{A}\int_0^\infty \int_0^a |T_{t,\alpha}f(s)-f(s)|^p t^{2\alpha+1}dt s^{2\alpha+1}ds$$ $$= \frac{1}{A}\int_0^a \int_0^\infty|T_{t,\alpha}f(s)-f(s)|^ps^{2\alpha+1}ds t^{2\alpha+1}dt.$$
Thus
\begin{equation}\label{mff}\|M_af-f\|_{p,(\alpha)}^p\le \sup_{0\le t\le a}\int_0^\infty|T_{t,\alpha}f(s)-f(s)|^ps^{2\alpha+1}ds.\end{equation}
Similarly

$$\|M_af-M_af_j\|_{p,(\alpha)}^p \le \left(\int_0^R|M_af(s)-M_af_j(s)|^p s^{2\alpha+1}ds\right)^{\frac{1}{p}}$$ $$+\left(\int_R^\infty |M_af(s)-M_af_j(s)|^p s^{2\alpha+1}ds\right)^{\frac{1}{p}}=I+II.$$
 As $II\le \left(\int_R^\infty |M_af(s)|^p s^{2\alpha+1}ds\right)^{\frac{1}{p}}+\left(\int_R^\infty |M_af_j(s)|^p s^{2\alpha+1}ds\right)^{\frac{1}{p}}$, it is enough to investigate $M_ag$, where $g\in K$.
$$\left(\int_R^\infty |M_ag(s)|^p s^{2\alpha+1}ds\right)^{\frac{1}{p}}\le \|M_ag-g\|_{p,(\alpha)}^p+\left(\int_R^\infty |g(s)|^p s^{2\alpha+1}ds\right)^{\frac{1}{p}}.$$
$$I\le \sup_{0\le s \le R}|M_af(s)-M_af_i(s)|B(R),$$
where $B(R)=\left(\frac{R^{2\alpha+2}}{2\alpha+2}\right)^{\frac{1}{p}}$. Thus
\begin{equation}\label{maf}\|M_af-M_af_j\|_{p,(\alpha)}^p \le B(R)\sup_{0\le s \le R}|M_af(s)-M_af_i(s)|\end{equation} $$+ \sup_{g\in K}\left(\|M_ag-g\|_{p,(\alpha)}^p+\left(\int_R^\infty |g(s)|^p s^{2\alpha+1}ds\right)^{\frac{1}{p}}\right).$$

Let $\varepsilon>0$ be arbitrary and let us choose by \eqref{c2} $R$ so large that $\forall f\in K$ $\left(\int_R^\infty|f(x)|^px^{2\alpha+1}dx\right)^{\frac{1}{p}}<\frac{\varepsilon}{5}$ and by \eqref{c3} let us choose $a$ so small such that for all $ 0\le h \le a$ and for all $f\in K$, \ws  $t\in[0,R] \ws \left(\int_0^\infty|T_{t+h,\alpha}f(s)-T_{t,\alpha}f(s)|^p s^{2\alpha+1}ds\right)^{\frac{1}{p}}<\frac{\varepsilon}{5}$. Now let us take an $\frac{\varepsilon}{5B(R)}$-net, $N:=N_{\frac{\varepsilon}{5B(R)}} \subset F_{a,R}$. By this choice for all $g\in K$
$$\|g-M_ag\|_{p,(\alpha)}\le \frac{\varepsilon}{5},$$
cf. \eqref{mff}. Moreover for an arbitrary $f\in K$ there is a suitable $f_j$ with $M_af_j \in N$, such that
$$\|M_af-M_af_j\|_{p,(\alpha)}\le B(R)\frac{\varepsilon}{5B(R)}+\frac{\varepsilon}{5}+\frac{\varepsilon}{5},$$
cf. \eqref{maf}. Finally we finish the proof with
$$\|f-f_j\|_{p,(\alpha)}\le \|f-M_af\|_{p,(\alpha)}+\|M_af-M_af_j\|_{p,(\alpha)}+\|M_af_j-f_j\|_{p,(\alpha)}\le \varepsilon.$$

\medskip

Now we define the Bessel-Fourier or Hankel transformation to  facilitate formulating a Pego-type theorem.  Let us denote the Hankel trasform of a function $f$ by
$$H_\alpha(f)(y)=\hat{f}(y):=c_\alpha\int_0^\infty f(x)j_\alpha(xy)x^{2\alpha+1}dx.$$
The inverse transformation is
$$H_\alpha^{[-1]}(\hat{f})(x)=c_\alpha\int_0^\infty \hat{f}(y)j_\alpha(xy)x^{2\alpha+1}dy,$$
if exists.
If $K$ is a set of certain functions, we denote by $\hat{K}$ the set of the Hankel transforms of the functions in  $K$.\\
Below we need that Hankel transformation fulfils the Hausdorff-Young inequality, that is
\begin{equation}\label{hy}  \|H_{\alpha} (f)\|_{p',(\alpha)} \le C_p \|f\|_{p,(\alpha)}, \ws \ws \ws  1\le p \le 2,\end{equation}
see \cite{gj}. We also use the next property of the Hankel transform.
\begin{equation}\label{ht} H_{\alpha} (T_{t,\alpha}f)(y)=j_\alpha(ty)\hat{f}(y); \ws \ws \ws T_{t,\alpha}\hat{f}(y)=H_\alpha(j_\alpha(t\cdot)f(\cdot))(y), \end{equation}
see \cite{gj}. Furthermore
\begin{equation}\label{hinf}\lim_{y\to\infty}\hat{f}(y)=0,\end{equation}
see \cite{p}. Defining the convolution by
$$ (f*g)(x)=\int_0^\infty T_{x,\alpha}f(t)g(t)t^{2\alpha+1}dt,$$
we have
\begin{equation}\label{hkonv} H_\alpha(f*g)=\hat{f}\hat{g},\end{equation}
see \cite{gj}.

\medskip

\begin{theorem}\label{He} $1\le p \le 2$. Let $K\subset L_{p,\alpha}$. If $K$ is  bounded and satisfies ${\bf P_A}$ in $L_{p,\alpha}$, $\hat{K}$ satisfies ${\bf P_B}$ in $L_{p',\alpha}$, and if $K$ satisfies ${\bf P_B}$ in $L_{p,\alpha}$, $\hat{K}$ satisfies ${\bf P_A}$ in $L_{p',\alpha}$. \end{theorem}

\proof
Let $K$ be bounded by $M$ and $f \in K$. Then by \eqref{ht} and then \eqref{hy}
$$\|T_{t+h,\alpha}\hat{f}-T_{t,\alpha}\hat{f}\|_{p',(\alpha)}=\|H_\alpha\left((j_\alpha((t+h)(\cdot))-j_\alpha(t(\cdot)))f(\cdot)\right)\|_{p',(\alpha)}$$ $$\le C_p \|(j_\alpha((t+h)s)-j_\alpha(ts))f(s)\|_{p,(\alpha)}\le C_p\left(\int_0^R \left|(j_\alpha((t+h)s)-j_\alpha(ts))f(s)\right|^ps^{2\alpha+1}ds\right)^{\frac{1}{p}}$$ $$+C_p\left(\int_R^\infty \left|(j_\alpha((t+h)s)-j_\alpha(ts))f(s)\right|^ps^{2\alpha+1}ds\right)^{\frac{1}{p}}=C_p(I+II).$$
Let $\varepsilon$ be arbitrary and we choose $R=R(\varepsilon)$ to $\frac{\varepsilon}{4}$ according to \eqref{c2}. By \eqref{j1}
$$II\le 2 \left(\int_R^\infty \left|f(s)\right|^ps^{2\alpha+1}ds\right)^{\frac{1}{p}}\le \frac{\varepsilon}{2}.$$
To estimate $I$, let $M_0>0$ be arbitrary and assume that $0<t<M_0$. Let us choose $\delta=\delta(R(\varepsilon),M_0,\varepsilon)<1$ such that
$$\delta \le \frac{(\alpha+1)\varepsilon}{(M_0+1)R^2M}.$$
Then in view of \eqref{j2} and \eqref{j1} for all $0\le h \le \delta$, by choice of $\delta$
$$I \le (M_0+1)R\frac{1}{2(\alpha+1)}\|j_{\alpha+1}\|_\infty \delta R \|f\|_{p,(\alpha)} \le \frac{\varepsilon}{2}.$$

To prove the second statement, for an arbitrary fixed $\delta>0$ we choose nonnegative functions $g_\delta$ supported on $(0,\delta)$ such that $$\int_0^\infty g_\delta(x)x^{2\alpha+1}dx=1.$$
Then in view of \eqref{hinf} we choose $R=R(\delta)$ such that $|H_\alpha(g_\delta)(y)|\le \frac{1}{2}$ if $y>R$. Thus, by \eqref{hkonv} and \eqref{hy}
$$\left(\int_R^\infty|\hat{f}(y)|^{p'}y^{2\alpha+1}dy\right)^{\frac{1}{p'}}\le 2 \left(\int_R^\infty|H_\alpha(f)(y)(1-H_\alpha(g_\delta)(y))|^{p'}y^{2\alpha+1}dy\right)^{\frac{1}{p'}}$$ $$\le C(p) \|f-f*g_\delta\|_{p,(\alpha)}=C(p)\left(\int_0^\infty\left|\int_0^\infty(f(x)-T_{x,\alpha}f(t))g_\delta(t)t^{2\alpha+1}dt\right|^px^{2\alpha+1}dx\right)^{\frac{1}{p}}$$ $$\le C(p)\int_0^\infty\left(\int_0^\infty|f(x)-T_{t,\alpha}f(x)|^px^{2\alpha+1}dx\right)^{\frac{1}{p}}g_\delta(t)t^{2\alpha+1}dt$$ $$\le C(p)\sup_{0\le t \le \delta}\left(\int_0^\infty|f(x)-T_{t,\alpha}f(x)|^px^{2\alpha+1}dx\right)^{\frac{1}{p}},$$
where in the last but one (Minkowski) inequality the symmetry of translation is considered. Thus the proof can be finished by choosing an appropriate $\delta$ to the arbitrary $\varepsilon$ in view of \eqref{c3}, and $g_\delta$ and $R(\delta)$ to $\delta$ as above.

\medskip

\begin{cor} $\mathrm{(I)}$ A bounded subset $K\subset L_{2,\alpha}$ is precompact if and only if for all $M_0>0$, $t\in [0, M_0]$\\
$\|T_{t+h,\alpha}f-T_{t,\alpha}f\|_{2,(\alpha)}\to 0$, and $\|T_{t+h,\alpha}\hat{f}-T_{t,\alpha}\hat{f}\|_{2,(\alpha)}\to 0$ uniformly in $K$.

\medskip

\noindent $\mathrm{(II)}$  A bounded subset $K\subset L_{2,\alpha}$ is precompact if and only if\\ $\left(\int_R^\infty|f(x)|^2x^{2\alpha+1}dx\right)^{\frac{1}{2}}\to 0$, and $\left(\int_R^\infty|\hat{f}(y)|^2y^{2\alpha+1}dy\right)^{\frac{1}{2}}\to 0$ uniformly in $K$.
\end{cor}

\medskip

\begin{cor} Let $u,v : \mathbb{R}_+ \to \mathbb{C}$, $u, v \in L_\infty$ such that $\lim_{x \to \infty}u(x), v(x)=0$, and assume that $u$ is continuous. Then the operator $L:= uH_\alpha v$, $f \mapsto uH_\alpha(vf)$ is compact on $L_{2,(\alpha)}$.
\end{cor}

\proof
Let $K\subset L_{2,\alpha}$ be a subset bounded by $M$. By the assumption, for an $\varepsilon>0$ there is an $R>0$ such that $|v(x)|<\varepsilon$ if $x>R$. Thus $\left(\int_R^\infty |(vf)(x)|^2x^{2\alpha+1}dx\right)^{\frac{1}{2}}\le \varepsilon M$, that is $vK$ is bounded and equivanishing. So by Theorem \ref{He} $H_\alpha(vK)$ is $L_{2,\alpha}$-equicontinuous, and by assumption $LK=uH_\alpha(vK)$ is also $L_{2,\alpha}$-equicontinuous and equivanishing.

\medskip

\remark
(1) In the proof of the first statement of Theorem \ref{He} we can also use the asymptotics of Bessel functions $|j_{\alpha+1}(u)|\le C u^{-\alpha+\frac{3}{2}}$. Then a larger $\delta$ can be allowed.\\
(2) Comparing \cite{hohm} to these special translations, the omission of the criterion of boundedness  needs further investigations in Bessel and Laguerre cases.

\medskip

\noindent \small{Department of Analysis, \newline
Budapest University of Technology and Economics}\newline

\small{ g.horvath.agota@renyi.hu}

\medskip

\begin{thebibliography}{99}

 \bibitem {a} R. A. Adams, {\it Sobolev spaces} Academic Press, New York (1975)

\bibitem {adh} V. Arestov, M. Deikalova, \'A. Horv\'ath, On Nikol'skii type inequality between the uniform norm and the integral q-norm with Laguerre weight of algebraic polynomials on the half-line, {\it J. of Approx. Theory}, {\bf 222} (2017) 40-54.

\bibitem {abdh} V. Arestov, A. Babenko, M. Deikalova, \'A. Horv\'ath, Nikol'skii inequality between the uniform norm and the integral norm with Bessel weight for entire functions of exponential type ont he half-line, {\it Analysis Math.}, {\bf 44} (2018) 21-42.

\bibitem {ag} R. Askey, G. Gasper, Convolution structures for Laguerre polynomials, {\it  J. Analyse Math.} {\bf 31} (1977) 48-68.

\bibitem {au} I. Aydin, C. Unal, The Kolmogorov-Riesz theorem and some compactness criterions of bounded subsets in weighted variable exponent amalgam and Sobolev spaces, {\it Collectanea Mathematica} {\bf 71} (2020) 349-367.

\bibitem {be} H. Bateman, A. Erd\'elyi, {\it  Higher Transcendental Functions, II}, McGraw-Hill (New York 1953).

\bibitem {bs} B. L. J. Braaksma and H. S. V. de Snoo, Generalized translation operators associated with a singular differential operator,{\it Proc. Conf. Theory of ordin. and part. diff. equat. in Dundee 1974} Springer Lecture Notes. {\bf 415} B. D. Sleeman, I. M.  Michael (Eds.).


\bibitem {gj} S. Ghobber, Ph. Jaming, Strong annihilating pairs for the Fourier-Bessel transform, {\it J. of Math. Analysis and Appl.} {\bf 377} (2011) 501-515.

\bibitem {dfg} M. D\"orfler, H. G. Feichtinger, K. Gr\"ochenig, Compactness criteria in function spaces, {\it Colloquium Mathematicae} {\bf 94} (2002) 37-50.

\bibitem{f} H. G. Feichtinger, A compactness criterion for translation invariant Banach spaces of functions, {\it Analysis Mathematica}, {\bf 8} (1982), 165-172.

\bibitem {g}  P. G\'orka, Pego theorem on locally compact Abelian goups, {\it Journal of Algebra and Its Applications}, {\bf 13}(2014) 1350143.

\bibitem {gk}  P. G\'orka, T. Kostrzewa, Pego everywhere, {\it Journal of Algebra and Its Applications}, {\bf 15}(2016) 1650074.

\bibitem {gma}  P. G\'orka, A. Macios, Almost everything you need to know about relatively compact sets in variable Lebesgue spaces, {\it Journal of Functional Analysis} {\bf 269} (2015)1925-1949.

\bibitem {gr} P. G\'orka, H. Rafeiro, From Arzel\'a-Ascoli to Riesz-Kolmogorov,  {\it Nonlinear Analysis} {\bf 144} (2016) 23-31.

\bibitem {gm} E. G\"orlich, C. Markett, A convolution structure for Laguerre series,  {\it Indag. Math.}, {\bf 44} (1982) 161-171.

\bibitem {gm1} E. G\"orlich, C. Markett, On approximation by Cesaro means of the Laguerre expansion and best approximation, {\it Resultate der Mathematik} {\bf 2} (1979) 124-150.

\bibitem {hoh} H. Hanche-Olsen, H. Holden, The Kolmogorov-Riesz compactness theorem, {\it Expositiones Mathematicae} {\bf 28} (2010) 385-94.

\bibitem {hohm} H. Hanche-Olsen, H. Holden, E. Malinnikova, An improvement of the Kolmogorov-Riesz compactness theorem, {\it Expositiones Mathematicae} {\bf 37} (2019) 84-91.

\bibitem {k} M. Krukowski, Characterizing compact families via Laplace transform, {\it Annales Acad. Sci. Fennicae Math.} {\bf 45} (2020) 991-1002.

\bibitem {ka} Y. Kanjin, On algebras with convolution structures for Laguerre polynomials, {\it Tran. of the Amer. Math. Soc.}, {\bf 295} (1986) 783-794.

\bibitem{le} B. M. Levitan, Expansion in Fourier series and integrals with Bessel functions, {\it Uspekhi Mat. Nauk} (in Russian) {\bf 6} ( 1951) 102--143.

\bibitem{m} C. Markett, Nikol'skii-type inequalities for Laguerre and Hermite expansions, in Colloquia Mathematica Societatis J\'anos Bolyai {\bf 35} (1980) 811-834.

 \bibitem{pe} R. L. Pego, Compactness in $L^2$ and the Fourier transform, {\it Proc. of the Amer. Math. Soc} {\bf 95} (1985) 252-254.

\bibitem{p} S. S. Platonov, Bessel harmonic analysis and approximation of functions on the half-line, {\it Izv. RAN, Ser. Mat.},  {\bf 71} (2007), 149-196 (in Russian); translated in {\it Izv. Math.}, {\bf 71}:5 (2007), 1001-1048.

\bibitem{r} H. Rafeiro, S. Samko, Dominated compactness theorem in Banach function spaces and its applications, {\it Complex Anal. Oper. Theory} {\bf 2} (4) (2008) 669-681.

\bibitem{s} V. N. Sudakov, Criteria of compactness in function spaces, (In Russian) {\it Upsekhi Math. Nauk.} {\bf 12} (1957) 221-224.

\bibitem {sz} G. Szeg\H o, {\it Orthogonal Polynomials}, AMS.  AMS Coll.~Publ.~Vol.~XXXIII (New York, 1959).

\bibitem {w} G. N. Watson, Another note on Laguerre polynomials, {\it J. London Math. Soc.} {\bf 14} (1939) 19-22.

\end{thebibliography}
\end{document}